# Synopsis of Duality

David L. Farnsworth
School of Mathematical Sciences
Rochester Institute of Technology
Rochester, New York 14623 USA

1. *Introduction*

Transformation of mathematical relationships is a powerful tool. Laplace and Fourier transformations are prominent examples. Taking logarithms, which can change a multiplication calculation into an addition calculation, may be classified as a transformation. We consider a specific transformation that takes points to lines and lines to points in one plane to another.

Consider two Euclidean planes. Call one plane the *original plane* and the other the *dual plane*. The original plane has rectangular coordinates $x$ and $y$, and the dual plane has rectangular coordinates $u$ and $v$.

The *dual transformation* between the two planes is that such the coordinates of a point in the original plane give the coefficients of a line in the dual plane, and the coefficients of a line in the original plane give the coordinates of a point in the dual plane. The point $(x,y) = (a,b)$ has dual line $au + bv = 1$, and the line $cx + dy = 1$ has dual point $(u,v) = (c,d)$. Since the correspondence is reversible, the dual of the dual of a line or plane is the original line or plane [1, pp. 24–25], [2, p. 143]. This property is called *reflexivity*. Since lines containing the origin cannot be written $cx + dy = 1$, they are either excluded from an analysis or their transformation is defined separately as necessary. The dual of a line that contains the origin can be taken to be a point at infinity, but such ideas are not necessary here. This exclusion is not as large a barrier as it might appear initially.

The main goal is to examine how the dual transformation carries curves from one plane to the other. In the original plane, consider a smooth curve that has a tangent line at each point and no linear portions. At a point on the curve, obtain the tangent line, and find its dual point, which is a point on the curve that is dual to the curve in the original plane. In detail, consider a region in the original $x,y$-plane where the curve can be written $y = y(x)$ with the tangent line at the point $(a, y(a))$ being

$$y = y(a) + y'(a)(x - a)$$

or

$$(-y'(a)/(y(a) - ay'(a)))x + (1/(y(a) - ay'(a)))y = 1.$$

The corresponding dual point is

$$(u(a), v(a)) = (-y'(a)/(y(a) - ay'(a)), 1/(y(a) - ay'(a))). \qquad (1)$$

A pair of parametric equations for the curve that is dual to $y = y(x)$ is given by (1). The line that is dual to the point $(x,y) = (a, y(a))$ is

$$au + y(a)v = 1. \qquad (2)$$

Equation (2) is the tangent line to the dual curve at the dual point (1), since its slope $-a/y(a)$ is

$$\frac{dv}{du} = \frac{\frac{dv}{da}}{\frac{du}{da}} = \frac{\frac{ay''(a)}{(y(a) - ay'(a))^2}}{\frac{-y(a)y''(a)}{(y(a) - ay'(a))^2}} = -\frac{a}{y(a)}. \qquad (3)$$

The curves in the two planes are dual to each other. In this way, there is a dual relationship between tangent lines, as well. If the parameter *a* can be eliminated from (1), the result is the dual curve $v = v(u)$.

This mutual duality between curves takes advantage of the fact that tangent lines characterize or determine sufficiently smooth curves, which are the envelopes of their tangent lines. At an isolated point where there is a discontinuity in the slope, such as a corner point at the intersection of two linear portions of a curve, there is no single tangent line. (We are considering linear portions of a curve to have tangent lines that contain them.) All lines that locally intersect the curve in just the corner point and are on only one side of the curve are used. Those lines are called *supporting lines* or *lines of support* [2, p. 41], [3, p. 37]. Examples of these points are (–2.5,0) and (2.5,0) in Figure 1 and the eight vertices in Figure 4. Curves with linear portions and many less-smooth curves pose no great difficulties, as illustrated by Examples 1, 4, and 5 in Section 2.

It is fruitful to identify the two dual planes. The *u* and *v* coordinates are renamed *x* and *y*, and all constructions and analyses are performed in a single plane. Except where it is noted, that identification is assumed below. This has the advantage that geometric constructions in the two planes can be directly compared and melded. For example, Theorem 3 offers a method for locating a point that is dual to a line, where the construction is in a single plane.

The examples and theorems show how the reflective dual relationship between points and lines in one plane and lines and points in another plane is extended to curves. For sufficiently smooth curves without linear portions, there is a point-by-point duality between curves in the two planes, which is illustrated by Examples 2 and 3. Tangent lines at dual points correspond, as well.

The idea of the dual plane is fundamental to linear algebra. The dot or scalar product

$$ux + vy = (u\ v)\begin{pmatrix} x \\ y \end{pmatrix}$$

contains $(x,y)$ from the original plane and $(u,v)$ from the dual plane of linear transformations. Both planes are vector spaces [1, p. 23], [4, pp. 73–90].

Section 2 contains examples of pairs of dual curves. Section 3 has some useful facts about the dual transformation. The Legendre transformation, which is based upon different coefficients of lines, is the subject of Section 4.

2. *Examples*

The examples are selected for their diversity and being as comprehensive as possible. Examples 1, 4, and 5 contain corners and linear portions. Example 2 gives the dual parabola for the parabola $y = x^2/(4p)$, which is in a standard form. The curves in Examples 3 and 4 are examples from convex analysis [5], [6, pp. 21–24], [7]. The inner curve in Example 4 is the so-called unit circle of taxicab geometry [8]. One curve in Example 5 has two points that share a line of support, which causes the surprising form of its dual curve.

*Example* 1: In Figure 1, the outer piecewise-differentiable curve consists of portions of the parabolas

$$y = -x^2/4 + 25/16 \text{ in the upper half-plane}$$

and

$$y = x^2/4 - 25/16 \text{ in the lower half-plane}.$$

The inner curve is its dual curve, which is composed of
$$x = \pm 0.4 \text{ for } -0.32 \le y \le 0.32,$$
$$y = 0.32 + 0.8(0.16 - x^2)^{1/2},$$
and
$$y = -0.32 - 0.8(0.16 - x^2)^{1/2}.$$

For a sample calculation, consider a point $(x_1, y_1) = (x_1, -x_1^2/4 + 25/16)$ on the upper portion of the outer curve. The tangent line is
$$y = (-x_1/2)x + x_1^2/4 + 25/16$$
or
$$((8x_1)/(4x_1^2 + 25))x + 16/(4x_1^2 + 25))y = 1.$$
The coefficients of $x$ and $y$ in the last equation are the coordinates of the point that is dual to the tangent line at $(x_1, y_1)$. Eliminating $x_1$ from these representations of the coordinates gives
$$y = 0.32 + 0.8(0.16 - x^2)^{1/2},$$
which is the upper portion of the dual curve. For example, the tangent line to the outer curve at $(x_1, y_1) = (1, 21/16)$ is
$$y = (-1/2)x + 29/16$$
or
$$(8/29)x + (16/29)y = 1,$$
which gives the dual point $(8/29, 16/29)$ on the inner curve.

Consider the point $(2.5, 0)$ on the outer curve, which has no tangent line, but has supporting lines
$$y = m(x - 2.5) + 0$$
or
$$(1/2.5)x - 1/(2.5m))y = 1$$
for $m \le -1.25$ and $m \ge 1.25$. These lines do not intersect the outer curve otherwise and are outside the curve. From the last equation, all of the $x$-coordinates of the dual points are $1/2.5 = 0.4$. The range of the $y$-coordinates is obtained from $y = -1/(2.5m)$ for $m \le -1.25$, which gives $0 \le y \le 0.32$, and for $m \ge 1.25$, which gives $-0.32 \le y \le 0$. This demonstrates how the dual of a corner point of a curve is a line segment in the dual curve, rather than a single point, and the dual of a line segment of a curve is one point.

As a point "moves" in the clockwise direction along either the inner or outer curve, its dual point moves in the same direction on its dual curve with possible exceptions of "pauses" or "instantaneous runs" as shown by the points of discontinuity of the slope of the outer curve and the corresponding linear portions of the inner curve.

*Example* 2: Figure 2 contains the two parabolas
$$y = x^2/(4p) \text{ and } y = -px^2$$
for $p > 0$, which are dual to each other. The tangent line at $(x_1, y_1) = (x_1, x_1^2/(4p))$ is
$$y = ((x_1/(2p))x - x_1^2/(4p)$$
or
$$(2/x_1)x - ((4p)/x_1^2)y = 1$$
for $x_1 \ne 0$. The dual point to $(x_1, y_1)$ is $(2/x_1, -(4p)/x_1^2)$. These points lie on $y = -px^2$, since $-(4p)/x_1^2 = -p(2/x_1)^2$.

To illustrate the inverse dual relationship, the tangent line at $(x_2, y_2) = (x_2, -px_2^2)$ is

$$y = -2px_2x + px_2^2$$

or

$$(2/x_2)x + (1/(px_2^2))y = 1,$$

for $x_2 \neq 0$. Since $1/(px_2^2) = (2/x_2)^2/(4p)$, these points lie on $y = x^2/(4p)$.

The parabolas' vertices at the origin are excluded from the duality relationship, but their positions and roles are clear.

Unlike Example 1, which contained corner points and linear segments, this duality is point-to-point or one-to one. Points in the first quadrant on $y = x^2/(4p)$ have dual points in the fourth quadrant on $y = -px^2$. Points in the second quadrant on $y = x^2/(4p)$ have dual points in the third quadrant on $y = -px^2$. As $(x,y)$ on $y = x^2/(4p)$ "moves" away from the origin in either direction, its dual point moves toward the origin.

*Example* 3: Consider the closed curve $|x|^p + |y|^p = 1$ with $p > 1$. Using its symmetries about the two axes, take $x > 0$ and $y > 0$, so that the absolute-value symbols may be omitted. The tangent line at $(x_1,y_1) = (x_1,(1 - x_1^p)^{1/p})$ is

$$y = (1 - x_1^p)^{1/p} - (x_1/(1 - x_1^p)^{1/p})^{p-1}(x - x_1)$$

or

$$x_1^{p-1}x + (1 - x_1^p)^{(p-1)/p}y = 1.$$

The point that is dual to the tangent line and hence to $(x_1,y_1)$ is $(x_1^{p-1},(1 - x_1^p)^{1-1/p})$. Adding the $p/(p-1)^{th}$ power of each coordinate shows that

$$(x_1^{p-1})^{p/(p-1)} + ((1 - x_1^p)^{(p-1)/p})^{p/(p-1)} = x_1^p + (1 - x_1^p) = 1.$$

The dual curve in the first quadrant is $x^{p/(p-1)} + y^{p/(p-1)} = 1$. Setting the exponent $p/(p-1) = q$ gives $1/p + 1/q = 1$ and $q > 1$. The curve dual to $|x|^p + |y|^p = 1$ is $|x|^q + |y|^q = 1$. If $p = 2$, then $q = 2$ and the two curves are unit circles and are *self-dual*. The dual curves for $p = 4$ (outer curve) and $q = 4/3$ (inner curve) are in Figure 3.

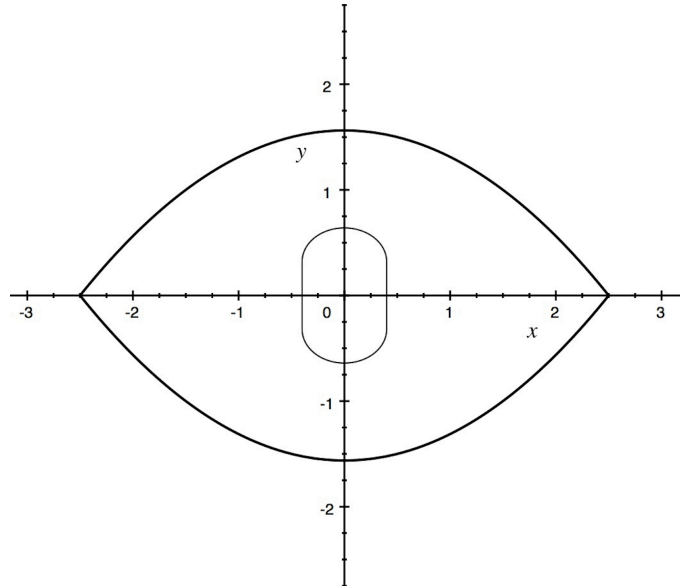

FIGURE 1: The two curves in Example 1.

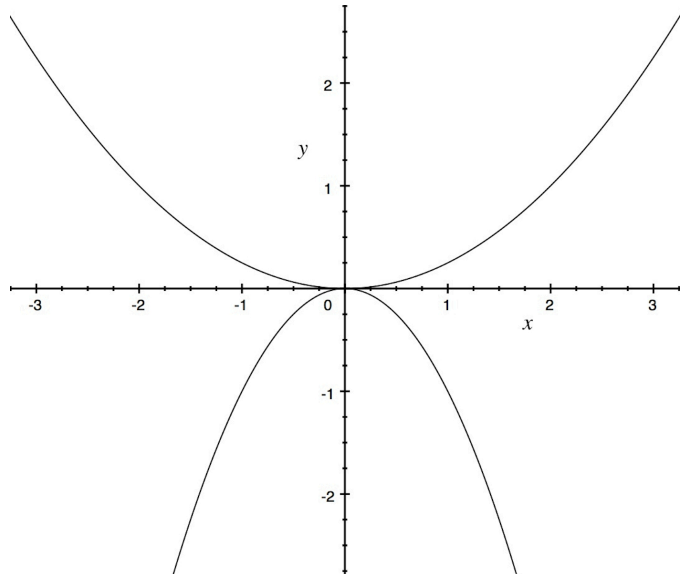

FIGURE 2: The two parabolas in Example 2 with $p = 1$.

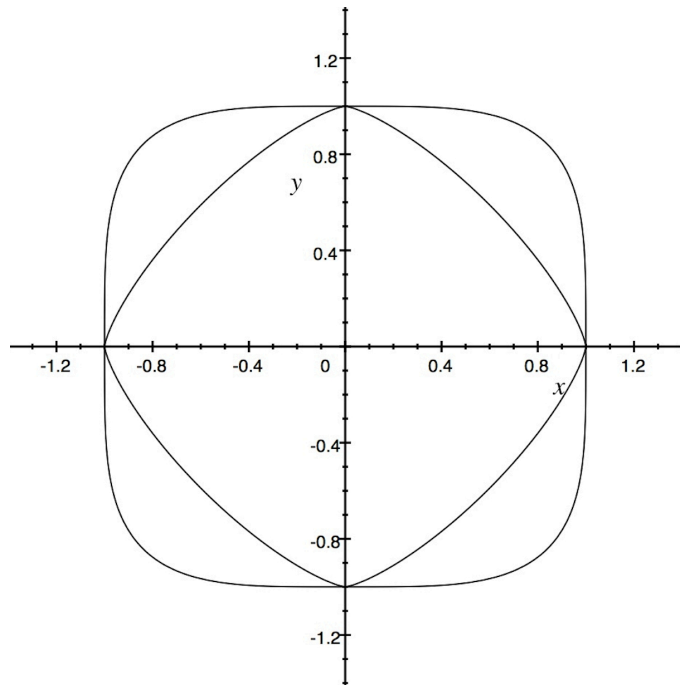

FIGURE 3: Curves in Example 3 with exponents 4 and 4/3.

*Example* 4: The curves are both piecewise linear in this example. Consider the closed, convex curve $|x| + |y| = 1$. Its dual curve is the square with vertices $(1,1)$, $(1,-1)$, $(-1,-1)$, and $(-1,1)$. The curves are in Figure 4. Using the symmetries about the two axes, take $x \geq 0$ and $y \geq 0$ as necessary, so that the absolute-value symbols are not required.

For the inner curve, all points on the side $x + y = 1$ have that side as their tangent line and have dual point (1,1), which is a vertex of the outer curve. For the inner curve's vertex (1,0), the supporting lines are $y = m(x - 1) + 0$ or $(-1)x + (-1/m)y = 1$ for $m \geq 1$ and $m \leq -1$, so that the coefficients are $-1$ and $-1/m$ with $-1 \leq -1/m \leq 1$. The dual of vertex (1,0) is the line segment $x = 1$ and $-1 \leq y \leq 1$, which is a side of the outer curve.

Conversely, for the outer curve, all points on its right-hand side $x = 1$ or $x + 0y = 1$ with $-1 \leq y \leq 1$ have dual point (1,0) on the inner curve, and all points on the side $y = 1$ or $0x + 1y = 1$ with $-1 \leq x \leq 1$ have dual point (0,1). The vertex (1,1) has supporting lines $y = m(x - 1) + 1$ or $(-m/(1 - m))x + (1/(1 - m))y = 1$ for $m \leq 0$. Hence, its dual on the inner curve is the portion of the line segment $x + y = 1$ with $0 \leq x \leq 1$.

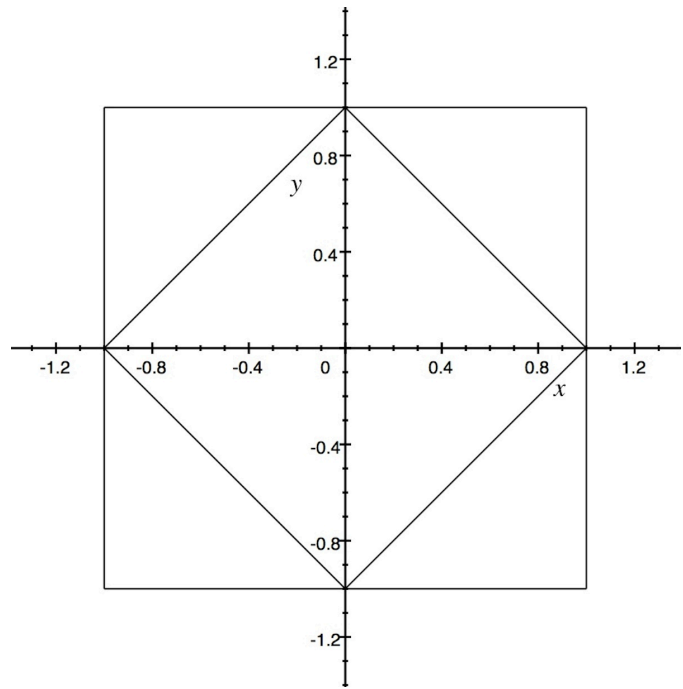

FIGURE 4: The two piecewise-linear curves in Example 4.

*Example* 5: This example illustrates that a dual curve exists for a curve when two points possess the same tangent line or supporting line. Computational details are omitted, since they are of the same kind as in the other examples. The curve in Figure 5 is

$$y = y(x) = \begin{cases} \sqrt{1 - x^2} & \text{for } -1 \leq x \leq 0 \\ -\dfrac{1}{3}x + 1 & \text{for } 0 \leq x \leq 1 \\ \dfrac{1}{3}x + \dfrac{1}{3} & \text{for } 1 \leq x \leq 2 \\ \sqrt{1 - (x - 2)^2} & \text{for } 2 \leq x \leq 3 \end{cases}$$

The two local maximum points at $x = 0$ and $x = 2$ have the same supporting line $y = 1$.

Its dual curve, which is displayed in Figure 6, is

$$y = y(x) = \begin{cases} \sqrt{1-x^2} & \text{for } -1 \le x \le 0 \\ 1 & \text{for } 0 \le x \le 1/3 \\ -\dfrac{3}{2}x + \dfrac{3}{2} & \text{for } -1 \le x \le 1/3 \\ 1 - 2x & \text{for } -1 \le x \le 0 \\ \sqrt{1 - 4x + 3x^2} & \text{for } 0 \le x \le 1/3 \end{cases}$$

As a point "moves" left to right in Figure 5, the dual point in Figure 6 follows the path *ABCDBE*. Point *B*, where the path crosses itself, has coordinates (0,1) from the shared supporting line $0x + 1y = 1$.

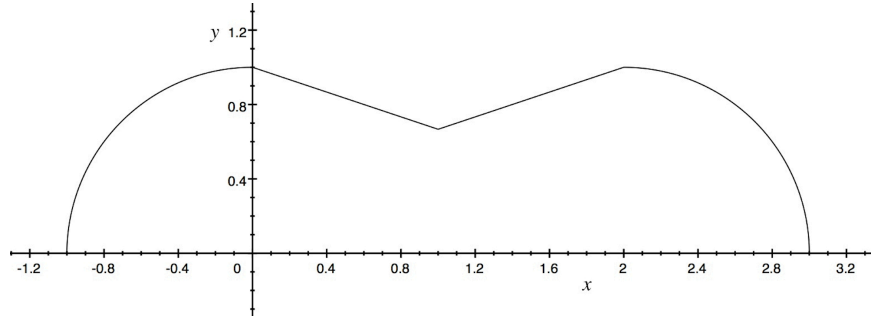

FIGURE 5: The original curve in Example 5.

3. *Geometry of duality*

This section presents selected facts about duality. It is understood that the dual transformation is not to be taken of lines through the origin.

*Theorem* 1: The dual lines of two points on line *L* intersect in a point, which is the dual of *L*. Two intersecting lines' dual points are on a line that is the dual of the two lines' point of intersection.

*Proof*: For any finite slope *m*, including $m = 0$, and non-zero *y*-intercept *b*, let *L* have the equation $y = mx + b$ or $(-m/b)x + (1/b)y = 1$. Two points on *L* are $P_1(x_1, mx_1 + b)$ and $P_2(x_2, mx_2 + b)$. The lines that are dual to $P_1$ and $P_2$ are $x_1 x + (mx_1 + b)y = 1$ and $x_2 x + (mx_2 + b)y = 1$, respectively, and both contain the point $(-m/b, 1/b)$, which is the dual of *L*. If *L* is the vertical line $x = x_1 \ne 0$ and the two points are $P_1(x_1, y_1)$ and $P_2(x_1, y_2)$, then the duals of $P_1$ and $P_2$ are $x_1 x + y_1 y = 1$ and $x_1 x + y_2 y = 1$, and both contain $(1/x_1, 0)$, which is the dual of *L*. ∎

*Theorem* 2: The dual points of two parallel lines with slope $m \ne 0$ are on $y = (-1/m)x$. For two parallel vertical lines, the dual points are on $y = 0$. For two parallel horizontal lines, the dual points are on $x = 0$.

*Proof*: If lines $L_1$ and $L_2$ have equations $y = mx + b_1$ and $y = mx + b_2$, that is, $(-m/b_1)x + (1/b_1)y = 1$ and $(-m/b_2)x + (1/b_2)y = 1$, with $m \ne 0$ and $b_1 b_2 \ne 0$, then the dual points to $L_1$ and $L_2$ are $(-m/b_1, 1/b_1)$ and $(-m/b_2, 1/b_2)$, which are on $y = (-1/m)x$. If the lines are the vertical lines $x = x_1$ and $x = x_2$ with $x_1 x_2 \ne 0$, then their dual points are $(1/x_1, 0)$ and $(1/x_2, 0)$, which lie on $y = 0$. If the lines are the horizontal lines $y = y_1$ and $y = y_2$ with $y_1 y_2 \ne 0$, then their dual points are $(0, 1/y_1)$ and $(0, 1/y_2)$, which lie on $x = 0$. ∎

In the proof of Theorem 2 for the first case with slope $m \neq 0$, if $L_1$ and $L_2$ are considered to meet at a point at infinity $P(\infty, m\infty)$ and "$\infty$" is treated algebraically as a number, then the dual of $P$ is $\infty x + m\infty y = 1$ or $y = (-1/m)x + 1/\infty = (-1/m)x$, which is the line that contains the dual points of the original parallel lines. This may be considered the limiting case of the second statement in Theorem 1. In Theorem 2, taking the limits of the $m \neq 0$ case as $m \to \infty$ and $m \to 0$ produce the outcomes for the vertical and for the horizontal parallel lines.

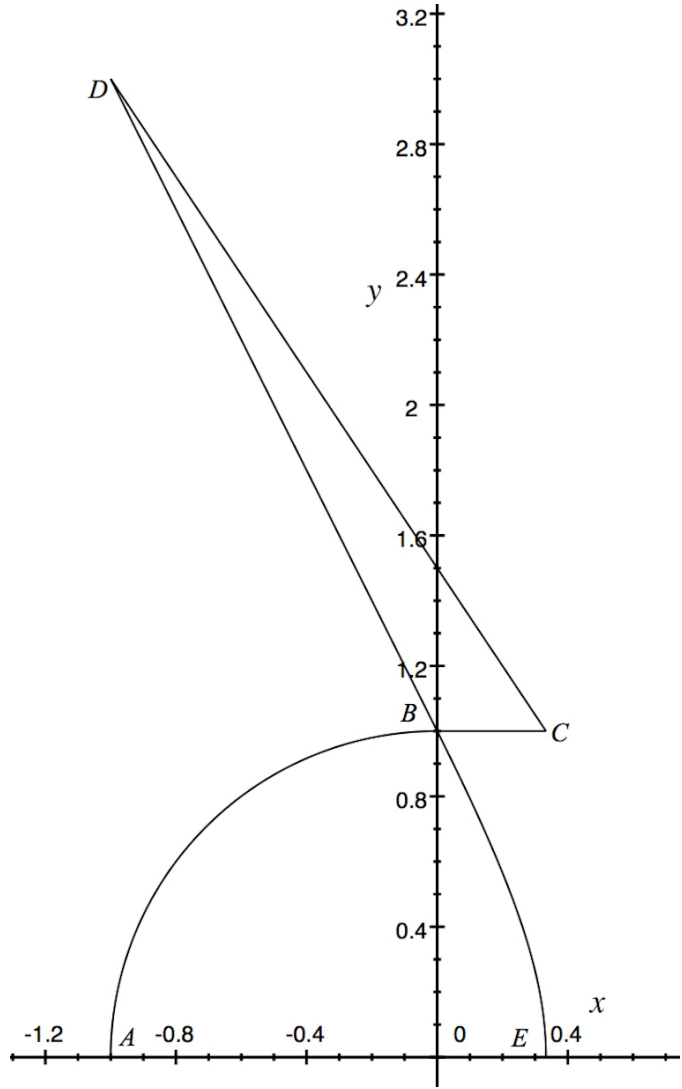

FIGURE 6: The dual curve in Example 5.

*Theorem* 3: The point $P$ that is dual to line $L_1$ is obtained by the following construction. Erect the line $L_2$ through the origin and perpendicular to $L_1$. Find the point $P_1$ at the intersection of $L_1$ and $L_2$. The point of inversion of $P_1$ with respect to the unit circle $x^2 + y^2 = 1$ is the dual point $P$ of $L_1$. Essentially reversing the steps gives the line $L_1$ that is dual to the point $P$.

*Proof*: The equation of $L_1$ may be written as $y = mx + b$, that is, $(-m/b)x + (1/b)y = 1$, whose dual is $P(-m/b, 1/b)$. See Figure 7. Line $L_2$ has equation $y = (-1/m)x$. The point of intersection of lines $L_1$ and $L_2$ is $(-bm/(m^2 + 1), b/(m^2 + 1))$. The squared distance between the origin $O$ and $P$ is
$$D^2(O,P) = m^2/b^2 + 1/b^2.$$
The squared distance between $O$ and $P_1$ is
$$D^2(O,P_1) = b^2m^2/(m^2 + 1)^2 + b^2/(m^2 + 1)^2.$$
Since $D^2(O,P)\, D^2(O,P_1) = 1$, $P$ and $P_1$ are inverse with respect to the unit circle. The same argument holds for vertical and horizontal lines $L_1$ that do not pass through $O$. ∎

*Corollary* 1: If line $L$ is tangent to a curve, then the radial direction to the dual point on the dual curve is perpendicular to $L$. The curve that is obtained by rotating the dual curve 90° about the origin has the property that each radial direction to a point is parallel to the tangent direction at its dual point on the original curve and each tangent direction at a point is parallel to the radial direction at its dual point on the original curve.

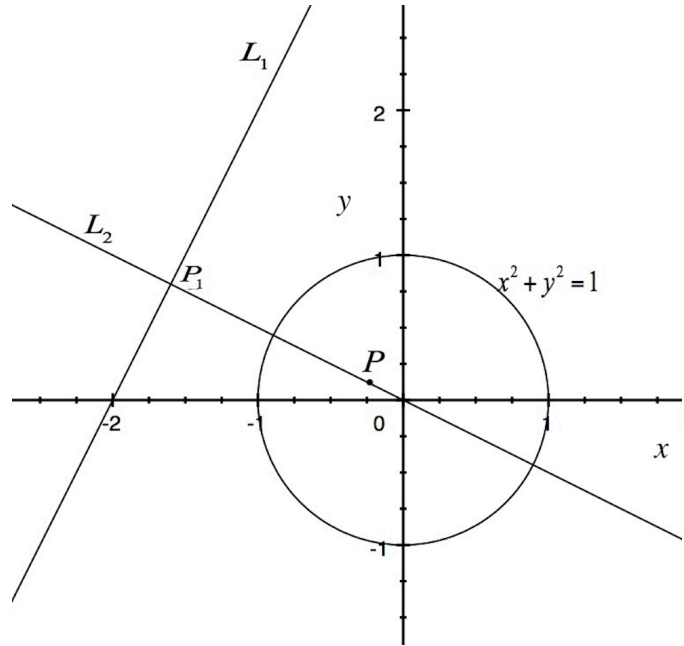

FIGURE 7: The construction in Theorem 3 for the point $P$ that is dual to $L_1$ and conversely for line $L_1$ that is dual to point $P$.

*Theorem* 4: Given lines $L_1$ and $L_2$ and their respective perpendicular lines $L_3$ and $L_4$ through the origin, the angle between $L_1$ and $L_2$ is equal to, the complement of, or the supplement of, the angle between $L_3$ and $L_4$.

*Proof*: Let $y = m_1x + b_1$ and $y = m_2x + b_2$ be the equations of $L_1$ and $L_2$, and let $y = (-1/m_1)x$ and $y = (-1/m_2)x$ be the equations of $L_3$ and $L_4$. The tangent of the angle between $L_1$ and $L_2$ is
$$(m_2 - m_1)/(1 + m_1m_2).$$
The tangent of the angle between $L_3$ and $L_4$ is
$$(1/m_1 - 1/m_2)/(1 + (1/(m_1m_2))) = (m_2 - m_1)/(1 + m_1m_2),$$

also. Figure 8 contains an example illustrating that the theorem is true if either $L_1$ or $L_2$ is horizontal or vertical. Alternatively, more geometrical arguments may be used. For example, consider the configuration that is displayed in Figure 9. In the right triangle *PQR*, the angle at *P* is $\theta$ and its complement at *R* is $\pi/2 - \theta$. Then, in the right triangle *RSO*, the angle at *R* is $\pi/2 - \theta$ and its complement at *O* is $\theta$. Figure 9 can be used to illustrate that line $L_1$ or $L_2$ can contain the origin and the theorem still holds. If $L_1$ contains the origin, then nearly the same figure may be used with $L_3$ perpendicular to $L_1$ through *O*. All of these steps are reversible. ∎

Theorem 4 is relevant to duality, because, according to Theorem 3, $L_3$ and $L_4$ contain the dual points to $L_1$ and $L_2$.

Theorem 5 shows that eliminating the parameter, such as $x_1$ in Examples 1–3, need not be performed in an *ad hoc* fashion.

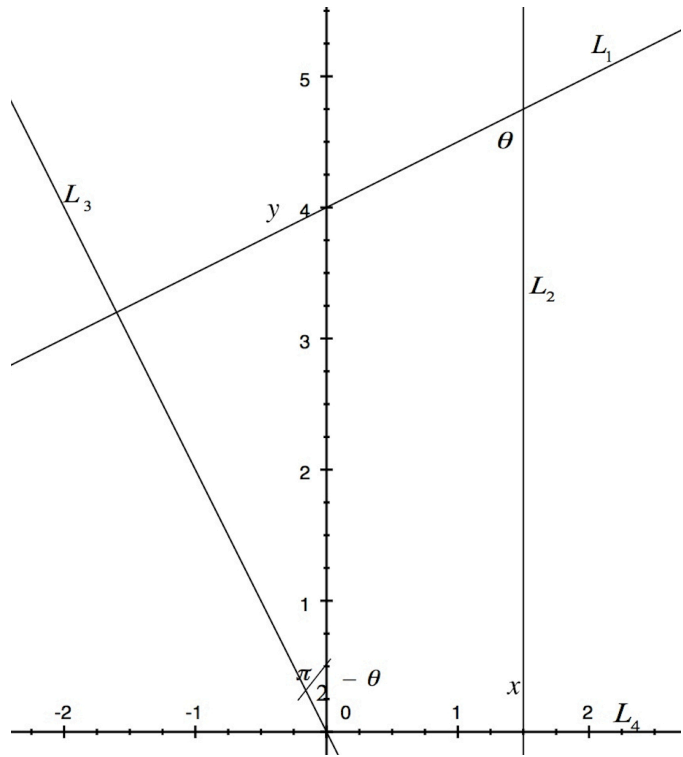

FIGURE 8: Example for Theorem 4 in which one line is vertical.

*Theorem* 5: For differentiable curves $y = y(x)$ for which $y'^{-1}$ is available, a formula for the dual curve is

$$y'^{-1}(-x/y)x + y(y'^{-1}(-x/y))y = 1. \qquad (4)$$

*Proof*: Take the derivative of (2) to obtain $x + y'(a)y = 0$, so that $a = y'^{-1}(-x/y)$ [9]. Substituting this into (2) gives (4). ∎

For readability, in Theorem 6 we revert to using the coordinates *u* and *v* in the dual plane.

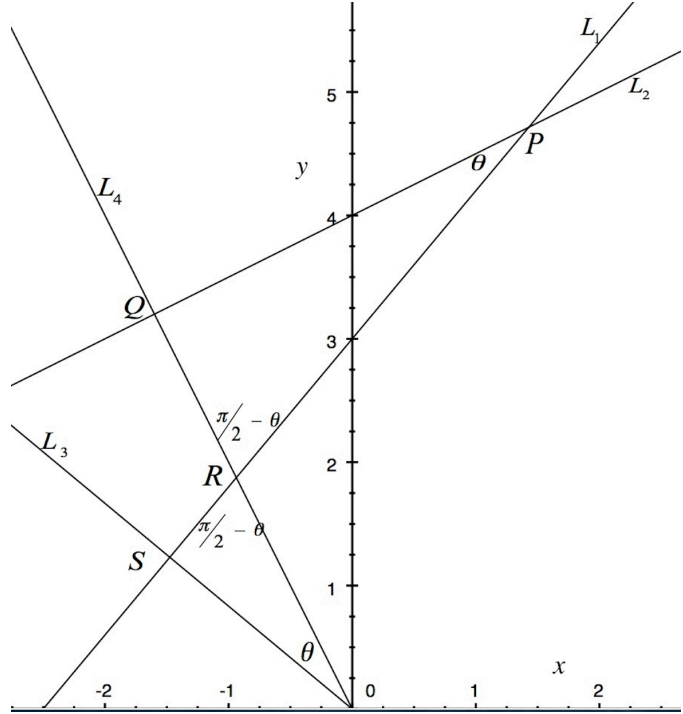

FIGURE 9: Two equal angles $\theta$ in Theorem 4.

*Theorem* 6: For twice differentiable functions $y = y(x)$ and $v = v(u)$ that represent dual curves, the second derivative of $v$ at the point which is dual to $(a, y(a))$ is

$$\frac{d^2v}{du^2} = \frac{(y(a) - ay'(a))^3}{(y(a))^3 y''(a)}. \tag{5}$$

*Proof*: Substituting the expressions for $u$ from (1) and for $dv/da$ from (2) into

$$\frac{d^2v}{du^2} = \frac{d}{da}\left(\frac{dv}{du}\right)\frac{da}{du} = \frac{\frac{d}{da}\left(\frac{dv}{du}\right)}{\frac{du}{da}}$$

yields (5). ∎

To illustrate Theorem 6, in Example 1, for the upper outer curve at $x = a$,
$$y(a) - ay'(a) > 0,\ y(a) > 0,\ \text{and}\ y''(a) < 0,$$
so (5) gives $d^2v/du^2 < 0$ at the dual point on the upper portion of the inner curve. In Example 2, for any real number $a$ and the parabola $y = x^2/(4p)$, $p > 0$,
$$y(a) - ay'(a) < 0,\ y(a) > 0,\ \text{and}\ y''(a) > 0,$$
so (5) gives $d^2v/du^2 < 0$ at the dual point.

*Theorem* 7: Expanding or contracting a curve by a factor $\lambda > 0$ expands or contracts the resultant dual curve by the factor $1/\lambda$.

*Proof* : For $(x_1, y(x_1))$ on the curve, the tangent line or a supporting line is $ux + vy = 1$ or

$$y = (-u/v)x + 1/v. \qquad (6)$$

Since the point is on the tangent or supporting line,
$$ux_1 + vy(x_1) = 1 \qquad (7)$$

For the corresponding point $(\lambda x_1, \lambda y(x_1))$ on the expanded or contracted curve, the tangent line or supporting line is $u_1 x + v_1 y = 1$ or
$$y = (-u_1/v_1)x + 1/v_1. \qquad (8)$$

Since the point is on the tangent or supporting line,
$$u_1(\lambda x_1) + v_1(\lambda y(x_1)) = 1. \qquad (9)$$

Since the tangent or supporting lines in (6) and (8) are parallel, their slopes are equal, implying that $-u_1/v_1 = -u/v$ or
$$v_1 = (u_1/u)v \qquad (10)$$

for $u \neq 0$. Substituting (10) into (9) gives
$$(x_1 u + y(x_1)v)\lambda u_1 = u.$$

Using (7) and (10) gives
$$u_1 = (1/\lambda)u \text{ and } v_1 = (1/\lambda)v. \blacksquare$$

Theorem 7 is reasonable from the point of view of the inversion construction in Theorem 3, where the coordinates of the point $P_1$ are multiplied by $\lambda$, and hence, the coordinates of the inverse point $P$ are divided by $\lambda$.

4. *The Legendre transformation*

Duality has been defined by using the coefficients of lines written as $ax + by = 1$. Other equations for lines yield other transformations. The relationship can be that the line $y = mx + b$ in one plane corresponds to the point $(m,-b)$ in the other plane. This transformation is called the *Legendre transformation* [10–12]. One advantage of the Legendre transformation is that lines through the origin with $m \neq 0$ are accommodated. The Legendre transformation is very useful in applications, especially in physics.

The Legendre transformation $\mathcal{L}\{y(x)\}$ of $y = y(x)$ is the function that is the negative of the tangent line's $y$-intercept as a function of its slope. Negation is introduced so that $\mathcal{L}$ is reflexive, that is, $\mathcal{L}\{\mathcal{L}\{y(x)\}\} = y(x)$. The present brief discussion is limited to sufficiently smooth functions with no linear portions and possessing derivatives with inverses. Each point $(a,y(a))$ corresponds to the pair $(m,t)$ by the Legendre transformation, where $y = mx - t$ is the tangent line at $(a,y(a))$. The Legendre transformation is $\mathcal{L}\{y(x)\} = t(m)$.

In order to express $t$ as a function of slope $m$ at the point $(a,y(a))$, write $m = y'(a)$ as $a = y'^{-1}(m)$. Then,
$$y(a) = y(y'^{-1}(m))$$
and
$$t(m) = ma - y(a) = my'^{-1}(m) - y(y'^{-1}(m)). \qquad (11)$$

By comparing the two representations $ux + vy = 1$ and $y = mx - t$ of a tangent line at a point, the simple change of variables $m = -u/v$, $t = -1/v$, that is, $u = m/t$, $v = -1/t$, between the dual curve and the Legendre transformation presents itself. Ignoring special cases where denominators are zero, techniques and algorithms can be borrowed between the two transformations. Sometimes, the Legendre transformation is called a dual transformation.

*Example* 3 *Revisited*: First, find $\mathcal{L}\{y(x)\}$. From $y = y(x) = (1-x^p)^{1/p}$ in the first quadrant,

$$y'(x) = -\frac{x^{p-1}}{(1-x^p)^{(1-p)/p}}.$$

Setting $y'(x) = m$, recalling that $x > 0$ and $m < 0$ in the first quadrant, and solving for $x$ gives

$$x = \frac{(-m)^{1/(p-1)}}{(1+(-m)^{p/(p-1)})^{1/p}} = y'^{-1}(m).$$

Substituting these into (11) yields the Legendre transformation

$$\mathcal{L}\{y(x)\} = t(m) = m\frac{(-m)^{1/(p-1)}}{(1+(-m)^{p/(p-1)})^{1/p}} - (1-(\frac{(-m)^{1/(p-1)}}{(1+(-m)^{p/(p-1)})^{1/p}})^p)^{1/p}$$

$$= -\frac{(-m)^{p/(p-1)}}{(1+(-m)^{p/(p-1)})^{1/p}} - (\frac{1+(-m)^{p/(p-1)}-(-m)^{p/(p-1)}}{1+(-m)^{p/(p-1)}})^{1/p}$$

$$= -\frac{1+(-m)^q}{(1+(-m)^q)^{1/p}} = -(1+(-m)^q)^{1-1/p} = -(1+(-m)^q)^{1/q}, \quad (12)$$

which is not part of the dual curve that is found in Example 3.
    Substituting $m = -x/y$, $t = -1/y$ into (12) gives
$$-1/y = -(1 + (-(-x/y))^q)^{1/q},$$

that is,
$$1/y = (1 + x^q/y^q)^{1/q}$$

or
$$x^q + y^q = 1,$$

which is the dual unit circle in the first quadrant.

*References*
1. P. R. Halmos, *Finite-dimensional vector spaces* (2nd edn.), Dover (2017).
2. S. R. Lay, *Convex sets and their applications*, Dover (2007).
3. R. V. Benson, *Euclidean geometry and convexity*, McGraw-Hill (1966).
4. D. C. Lay, *Linear algebra and its applications* (3rd edn.), Pearson/Addison-Wesley (2006).
5. V. Klee, What is a convex set?, *American Mathematical Monthly* **78** (1971) pp. 616–631.
6. A. C. Thompson, *Minkowski geometry*, Cambridge (1996).
7. J. B. Keller and R. Vakil, $\pi_p$, the value of $\pi$ in $\ell_p$, *American Mathematical Monthly* **116** (2009) pp. 931–935.
8. E. F. Krause, *Taxicab geometry: an adventure in non-Euclidean geometry*, Dover (1986).
9. W. R. Longley, Note on a theorem of envelopes, *Annals of Mathematics, Second Series*, **17** (1916), pp. 169-171.

10. S. Kennerly, A graphical derivation of the Legendre transform, 9 pages, available at www.physics.drexel.edu/~skennerly/maths/Legendre.pdf, assessed July 10, 2018.

11. C. E. Mungan, Legendre transforms for dummies, 11 pages, available at www.aapt.org/docdirectory/meetingpresentations/SM14/Mungan-Poster.pdf, assessed July 10, 2018.

12. R. K. P. Zia, E. P. Redish, and S. R. McKay, Making sense of the Legendre transform. *American Journal of Physics* **77** (2009) pp. 614–622.14